\documentclass[12pt]{article}
\usepackage{amsmath, amssymb, eucal, latexsym}
\usepackage[cp1251]{inputenc}
\usepackage[english]{babel}
\usepackage{graphicx}

\def\R{\Bbb R}

\def\a{\alpha}
\def\la{\lambda}
\def\D{\Delta}

\def\E{\mathsf {E}}

\def\Sym{\mathsf {Sym}}

\def\d{\delta}

\newtheorem{theorem}{Theorem}

\newtheorem{lemma}[theorem]{Lemma}
\newtheorem{cor}[theorem]{Corollary}
\newtheorem{note}[theorem]{Remark}
\newtheorem{definition}[theorem]{Definition}

\begin{document}

 \begin{center}
% \textbf{+ суммах множеств, and меіих малое произведение}
 \textbf{NEW RESULTS ON SUM--PRODUCTS IN $\R$}
 \end{center}

 \begin{center}
                                                         S. V. KONYAGIN, I. D. SHKREDOV
\footnote{
This work is supported by the Russian Science Foundation under a grant 14-50-00005.
%іабота вvполнена в іатематическом and нститут and м. T.L. Tтеклова іL=.
}\\

    \end{center}

\bigskip

\begin{center}
    Abstract.
\end{center}

{\it \small
    We improve a previous sum--products estimates in $\R$, namely, we obtain that
    $\max{\{|A+A|,|AA|\}}\gg |A|^{\frac{4}{3}+c}$, where $c$ any number less than $\frac{5}{9813}$.
    %Also
    New lower bounds for sums of sets with small
    the product set
    %difference set
    are found.
    Also we prove some pure energy sum--products results, improving a result of Balog and Wooley, in particular.
}

\bigskip
\section{Introduction}
\label{sec:introduction}
\bigskip

%SK
Let  $A,B\subset \R$ be finite sets.
Define the  \textit{sum set}, the \textit{product set} and \textit{quotient set} of $A$ and $B$ as
$$A+B:=\{a+b ~:~ a\in{A},\,b\in{B}\}\,,$$
$$AB:=\{ab ~:~ a\in{A},\,b\in{B}\}\,,$$
and
$$A/B:=\{a/b ~:~ a\in{A},\,b\in{B},\,b\neq0\}\,,$$
correspondingly.
The Erd\"{o}s--Szemer\'{e}di  conjecture \cite{ES} says that for any  $\epsilon>0$ one has
$$\max{\{|A+A|,|AA|\}}\gg{|A|^{2-\epsilon}} \,.$$
%SK +рубо говор , утверждаетс , что произвольное set веественнvх (или целvх) чисел не может одновременно
%SK and меть
Roughly speaking, it asserts that an arbitrary subset of real numbers (or integers)
cannot has good additive and multiplicative structure, simultaneously.
%At the moment the best result in this  direction is due to
Using some beautiful geometrical arguments Solymosi \cite{soly}, proved the following

\begin{theorem}
    Let  $A\subset \R$ be a set.
    % and $\tau>0$ --- веественное число.
    Then
\begin{equation}\label{f:Solymosi}
%SK    |A+A|^2 |AA| \,, \quad  |A+A|^2 |A/A| \ge \frac{|A|^4}{4 \lceil \log |A| \rceil} \,.
    |A+A|^2 |A/A|\,, \quad |A+A|^2 |AA|   \ge \frac{|A|^4}{4 \lceil \log |A| \rceil} \,.
\end{equation}
    In particular
\begin{equation}\label{f:Solymosi_max}
    \max{\{|A+A|,|AA|\}}\gg \frac{|A|^{4/3}}{\log^{1/3} |A|} \,.
\end{equation}
\label{t:Solymosi}
\end{theorem}

Here and below we suppose that $|A|\ge2$.

     It is easy to see that bound (\ref{f:Solymosi}) is tight up to logarithmic factors if
     the size of $A+A$ is small relatively  to $A$.
     We will write $a \lesssim b$ or $b \gtrsim a$ if $a = O(b \cdot \log^c |A|)$, $c>0$.
     If $a \lesssim b$ and $b \lesssim a$ then we write $a\sim b$.

\bigskip

In
%Theorem  \ref{t:Sol+}
paper  \cite{KS_smd} we improved bound (\ref{f:Solymosi_max}).

\begin{theorem}
    Let  $A\subset \R$ be a set.
    Then
$$\max{\{|A+A|,|AA|\}} \gtrsim |A|^{\frac{4}{3}+c'} \,,$$
where $c'=\frac{1}{20598}$.
% is any constant.
The same is true if one replace $AA$ by $A/A$.
\label{t:Sol_new}
\end{theorem}

The main result of the article is the following.

\begin{theorem}
    Let  $A\subset \R$ be a set.
    Then
$$\max{\{|A+A|,|AA|\}} \gtrsim |A|^{\frac{4}{3}+c} \,,$$
where $c=\frac{5}{9813}$.
% is any constant.
The same is true if one replace $AA$ by $A/A$.
\label{t:main2}
\end{theorem}

In paper \cite{KS_smd} the case of sets with small the product/quotient  sets was considered
(sharper bounds for {\it difference} of two sets, having small multiplicative doubling can be found in \cite{Sh_ineq}).

\begin{theorem}
    Let $A\subset \R$ be a finite set and $K\ge 1$ be a real number.
%SK    Suppose that $|AA| \le K|A|$ or $|A/A| \le K|A|$.
Suppose that $|A/A| \le K|A|$ or $|AA| \le K|A|$.
    Then
\begin{equation}\label{f:main_1}
    |A+A| \gtrsim |A|^{\frac{19}{12}} K^{-\frac{5}{6}}
    %\,.
\end{equation}
    and
\begin{equation}\label{f:main_2}
 |A+A| \gtrsim |A|^{\frac{49}{32}} K^{-\frac{19}{32}} \,.
\end{equation}
\label{t:main}
\end{theorem}

Inequality (\ref{f:main_2}) is better than (\ref{f:main_1}) for $K \gtrsim |A|^{\frac{5}{23}}$.

%Let us formulate the first result of the article
%(its  refined version is contained in  Theorem \ref{t:small_md} and Theorem \ref{t:small2_md} below).

%Theorem \ref{t:main} is stronger than  Theorem \ref{t:previous} and
%refines estimate (\ref{f:sol}) for $K \lesssim |A|^{1/3}$.

%Also,
We improve Theorem \ref{t:main} for some range of parameters in the case of small quotient set.

\begin{theorem}
    Let $A\subset \R$ be a finite set and $K\ge 1$ be a real number.
%SK    Suppose that $|AA| \le K|A|$ or $|A/A| \le K|A|$.
Suppose that $|A/A| \le K|A|$.
% or $|AA| \le K|A|$.
    Then
\begin{equation}\label{f:main1}
    |A+A| \gtrsim \max\{ |A|^{\frac{19}{12}} K^{-\frac{5}{6}}, |A|^{\frac{1313}{830}} K^{-\frac{336}{415}} \} \,.
    %\,.
\end{equation}
\label{t:main1}
\end{theorem}

One can check that lower bound (\ref{f:main1}) coincides with (\ref{f:main_1}) for $K \lesssim |A|^{\frac{5}{23}}$
and is better than both estimates (\ref{f:main_1}), (\ref{f:main_2}) for $|A|^{\frac{5}{23}} \lesssim  K \lesssim |A|^{\frac{673}{2867}}$. If $K \gtrsim |A|^{\frac{673}{2867}}$ then  (\ref{f:main_2}) gives better result.

Finally, in section \ref{sec:sp-energies} we
prove
sum--products results, which have deal just with the energies of sets but not with its sumsets or product sets.
Similar results in the direction were obtained in \cite{BW}, where the following Theorem was proved.

\begin{theorem}
    Let $A\subset \R$ be a finite set and $\d = 2/33$.
    Then there are two disjoint subsets $B$ and $C$ of $A$ such that $A = B\sqcup C$ and
$$
    \max\{ \E^{+} (B), \E^{\times} (C)\} \ll |A|^{3-\delta} (\log |A|)^{1-\delta}
$$
    and
$$
    \max\{ \E^{+} (B,C), \E^{\times} (B,C)\} \ll |A|^{3-\d/2} (\log |A|)^{(1-\d)/2} \,.
$$
\label{t:BW}
\end{theorem}

Also it was proved in \cite{BW} that one cannot take $\d$ 
%less 
greater 
than $2/3$.
Our method gives an
%small
improvement of Theorem \ref{t:BW}.
% of lower bound for $\d$.

\begin{theorem}
    Let $A\subset \R$ be a finite set and $\d = 1/5$.
    Then there are two disjoint subsets $B$ and $C$ of $A$ such that $A = B\sqcup C$ and
$$
    \max\{ \E^{+} (B), \E^{\times} (C)\} \lesssim  |A|^{3-\d} \,.
$$
\label{t:BW'}
\end{theorem}

In the proof of our results we use a combination of methods from \cite{soly}, \cite{SS1} and of course  \cite{KS_smd}.
 %in our arguments.
The main additional idea is to introduce some more flexible quantity $d_* (A)$ instead of quantity $d(A)$, see the definitions below. It allows us to avoid of using the Balog--Szemer\'{e}di--Gowers Theorem \cite{TV}.
%which
This usually provides better bounds and allows us, in addition, to obtain a series of
pure energy
results in section \ref{sec:sp-energies}.
We hope that our new quantity $d_* (A)$ will help in another problems of sum--products type.

\bigskip
\section{Definitions and preliminary results}
\label{sec:definitions}
\bigskip

%=апомним необходимvе определени .
The {\it additive energy $\E^{+} (A,B)$} between two sets $A$ and $B$ is the number of the solutions of the equation
(see \cite{TV})
$$
    \E^{+} (A,B) = |\{ a_1+b_1 = a_2+b_2 ~:~ a_1,a_2\in A\,, b_1,b_2\in B \}| \,.
$$
The {\it multiplicative energy $\E^{\times} (A,B)$} between two sets $A$ and $B$ is  the number of the solutions
of the equation (see \cite{TV})
$$
    \E^{\times} (A,B) = |\{ a_1 b_1 = a_2 b_2 ~:~ a_1,a_2\in A\,, b_1,b_2\in B \}| \,.
$$
In the case $A=B$ we write $\E^{+} (A)$ for $\E^{+} (A,A)$ and $\E^{\times} (A)$ for $\E^{\times} (A,A)$.
%SK -л  $\la\in A/A$ обозначим $A_\la = A\cap \la A$. -сно, что
Having $\la\in A/A$, we  put $A_\la = A\cap \la A$. Clearly, if $0\not\in A$ then
\begin{equation}\label{f:energy_basic}
    \E^{\times} (A) = \sum_{\la \in A/A} |A_\la|^2
    %\,.
\end{equation}
and, similarly, for the energy $\E^{+} (A)$.
%%+сли $P\subseteq A/A$, то пишем $\E^{\times}_P (A)$ for $\sum_{\la \in P} |A_\la|^2$.
Next, the Cauchy--Schwarz inequality implies for $0\not\in A$, $A_1\subset A$, $A_2\subset A$ that
\begin{equation}\label{f:energy_KB_subs}
%SK    \E^{\times} (A_1,A_2) |AA| \ge |A_1|^2|A_2|^2 \,, \quad \quad
    \E^{\times} (A_1,A_2) |A/A| \ge |A_1|^2|A_2|^2 \,, \quad \quad
\E^{\times} (A_1,A_2) |AA| \ge |A_1|^2|A_2|^2 \,.
\end{equation}
In particular
\begin{equation}\label{f:energy_KB}
%SK    \E^{\times} (A) |AA| \ge |A|^4 \,, \quad \quad \E^{\times} (A) |A/A| \ge |A|^4 \,.
    \E^{\times} (A) |A/A| \ge |A|^4 \,, \quad \quad \E^{\times} (A) |AA| \ge |A|^4 \,.
\end{equation}
Finally, we will use the following inequality.
\begin{lemma}\label{1/4ineq}
Let $A_1,\dots,A_n$ be finite subsets of $\R$. Then
$$
    \left(\E^{+} \left(\bigcup_{i=1}^n A_i\right)\right)^{1/4} \le \sum_{i=1}^n (\E^{+} (A_i))^{1/4}.
$$
Similarly, if $A_1,\dots,A_n$ are finite subsets of $\R\setminus\{0\}$, then
$$
    \left(\E^{\times} \left(\bigcup_{i=1}^n A_i\right)\right)^{1/4} \le \sum_{i=1}^n (\E^{\times} (A_i))^{1/4}.
$$
\end{lemma}

{\bf Proof.}
A similar result for subsets of finite abelian groups follows from inequality (4.18)
and Exercise 4.2.1 from \cite{TV}. Subsets of $\R$ can be reduced to subsets of finite
groups by Lemma~5.26 from \cite{TV}.
$\hfill\Box$

\bigskip

We need in several
%lemmas.
%results.
auxiliary statements.
The first  one is
%a consequence of
the Szemer\'edi--Trotter Theorem
\cite{sz-t}, see also \cite{TV}. We call a set $\mathcal{L}$ of continuous
plane curves a {\it pseudo-line system} if any two members of $\mathcal{L}$
share at most one point in common.
Define the {\it number of indices} $\mathcal{I} (\mathcal{P},\mathcal{L})$ between points and pseudo--lines  as
$\mathcal{I}(\mathcal{P},\mathcal{L})=|\{(p,l)\in \mathcal{P}\times \mathcal{L} : p\in l\}|$.

\begin{theorem}\label{t:SzT}
%(\cite{sz-t})
Let $\mathcal{P}$ be a set of points and let $\mathcal{L}$ be a pseudo-line system.
Then
$$\mathcal{I}(\mathcal{P},\mathcal{L}) \ll |\mathcal{P}|^{2/3}|\mathcal{L}|^{2/3}+|\mathcal{P}|+|\mathcal{L}|\,.$$
\end{theorem}

\bigskip

We need in a definition,  see \cite{s_sumsets}.

\begin{definition}
%[{\bf Shkredov \cite{s_sumsets}}]
A finite set $A \subset \R$ is said to be of {\it Szemer\'{e}di--Trotter type}
%SK  (abbreviated as {\it SzT--type}) if there exists a parameter $D(A)>0$ such that inequality
(abbreviated as {\it SzT--type}) with a parameter $D>0$ if the inequality
\begin{equation}\label{f:SzT-type}
\bigl|\bigl\{ s\in A-B ~\mid~ |A\cap (B+s)| \ge \tau \big\}\big|
    \le
%        \ll
%SK        \frac{D(A) |A| |B|^2}{\tau^{3}}\,,
\frac{D |A| |B|^2}{\tau^{3}}\,,
\end{equation}
holds
 for every finite set $B\subset \R$ and every real number $\tau \ge 1$.
\label{def:SzT-type}
\end{definition}

%SK So, $D(A)$ can be considered as the infimum of numbers such
 The quantity $D(A)$ can be considered as the infimum of numbers $D$ such that (\ref{f:SzT-type})
 %holds
 takes place
 for any $B$ and $\tau \ge 1$
 but, of course, the definition is applicable just for sets $A$ with  small quantity $D(A)$.

%Lemma \ref{l:d(A)} implies the following  result.
Any SzT--type set has small number of solutions of a wide class of linear equations, see e.g. Corollary 8 from \cite{KS_smd} (where nevertheless another quantity $D(A)$ was used)
%or
and
Lemma 7, 8 from \cite{s_sumsets}, say.

\begin{cor}
    Let  $A_1,A_2,A_3 \subset \R$ be any finite sets
    and $\a_1,\a_2,\a_3$ be arbitrary nonzero numbers.
    Then the number of the solutions of the equation
\begin{equation}\label{f:gen_sigma}
    \sigma (\a_1 A_1, \a_2 A_2, \a_3 A_3) :=
    |\{ \a_1 a_1 + \a_2 a_2 + \a_3 a_3 = 0 ~:~ a_1\in A_1\,, a_2\in A_2\,, a_3\in A_3 \}|
\end{equation}
    does not exceed
    $O(D^{1/3} (A_1) |A_1|^{1/3} |A_2|^{2/3} |A_3|^{2/3})$.
    Further, $\E^{+} (A_1,A_2) \ll D^{1/2} (A_1) |A_1| |A_2|^{3/2}$.
\label{c:gen_sigma}
\end{cor}

Also we need in a result from \cite{s_sumsets} on connection between the sumsets and $D(A)$ for SzT--type sets $A$.

\begin{theorem}
    Let $A$ has SzT type.
    Then
\begin{equation}\label{f:s_SzT}
    |A+A| \gtrsim |A|^{\frac{58}{37}} D(A)^{-\frac{21}{37}} \,.
\end{equation}
\label{t:s_SzT}
\end{theorem}

Now we can introduce a new characteristic of a set $A\subset \R$.
Put
$$
    \Sym^\times_{t} (Q,R) = \{ x ~:~ |Q \cap xR^{-1}| \ge t \} \,,
$$
and
\begin{equation}\label{f:d_r}
    d_* (A) = \min_{t>0}\, \min_{\emptyset \neq Q,R \subset \R \setminus \{0\}
                    % ~:~ A\subseteq \Sym^\times_{t} (Q,R)
                    } \,
        \frac{|Q|^2 |R|^2}{|A| t^3} \,,
\end{equation}
    where the second minimum in (\ref{f:d_r}) is taken over any $Q,R$ such that $A\subseteq \Sym^\times_{t} (Q,R)$
    and $\max\{ |Q|,|R| \} \ge |A|$.

\begin{lemma}
    Let $A\subset \R$ be a finite set.
    Then $A$ is of Szemer\'{e}di--Trotter type with $O(d_*(A))$.
\label{l:d_r}
\end{lemma}
{\bf Proof.}
    Let $R$, $Q$ be two sets and $t>0$ be a real number such that $A\subseteq \Sym^\times_{t} (Q,R)$.
    Without loosing of generality assume that $|Q| = \max\{ |Q|,|R| \} \ge |A|$.
    Let also
$$
    S_\tau := \{ s \in A-B ~:~  |A\cap (B+s)| \ge \tau \} \,.
$$
    Our task is to estimate the size of $S_\tau$.
    It is easy to see that the  bound
\begin{equation}\label{f:S_tau_RQ}
    |S_\tau| \ll \frac{|Q|^2 |R|^2 |B|^2}{t^3 \tau^3} \,.
\end{equation}
    is enough.
    We have
\begin{equation*}\label{tmp:02.12.2015_1}
    \tau |S_\tau| \le \sum_{s\in S_\tau} |A\cap (B+s)| = |\{ a-b = s ~:~ a\in A,\, b\in B,\, s\in S_\tau \}|
        := \sigma \,.
\end{equation*}
Because $A\subseteq \Sym^\times_{t} (Q,R)$, we obtain the following upper bound for the number of solutions $\sigma$
\begin{equation}\label{tmp:03.12.2015_1-}
    \sigma \le t^{-1} |\{ qr-b = s ~:~ q\in Q,\, r\in R,\, b\in B,\, s\in S_\tau \}|
\end{equation}
First of all let us prove a trivial estimate for the size of $S_\tau$.
Namely, dropping the condition $s\in S_\tau$ in (\ref{tmp:03.12.2015_1-}),
%summing (\ref{tmp:03.12.2015_1}) over all $s$,
we get
$$
    \tau |S_\tau| t \le |Q| |R| |B|
$$
and hence inequality (\ref{f:S_tau_RQ}) should be checked in the range
\begin{equation}\label{f:t_tau_range}
    t^2 \tau^2 \gg |Q| |R| |B|
\end{equation}
only because otherwise
$$
    |S_\tau| \le \frac{|Q||R||B|}{t\tau} \ll \frac{|Q|^2 |R|^2 |B|^2}{t^3 \tau^3} \,.
$$

Further, consider the family $\mathcal{L}$ of $|R||S_\tau|$ lines $l_{r,s} = \{ (x,y) ~:~ ry - x = s\}$, $r\in R$, $s\in S_\tau$
and the family of points $\mathcal{P} = Q \times B$.
Applying Theorem \ref{t:SzT} to
%the families,
the pair $(\mathcal{P}, \mathcal{L})$,
we get
\begin{equation}\label{tmp:03.12.2015_1}
    \sigma \le t^{-1} \mathcal{I} (\mathcal{P}, \mathcal{L})
        \ll
            t^{-1} \left( (|\mathcal{P}| |\mathcal{L}|)^{2/3} + |\mathcal{P}| + |\mathcal{L}| \right)
\end{equation}
If the first term in (\ref{tmp:03.12.2015_1}) dominates then we obtain (\ref{f:S_tau_RQ}).
Now suppose that required bound (\ref{f:S_tau_RQ}) does not hold.
Then if the second term in (\ref{tmp:03.12.2015_1})
%dominates,
is the largest one,
we obtain
$$
    \frac{|Q|^2 |R|^2 |B|^2}{t^2 \tau^2} \ll t \tau |S_\tau| \ll |\mathcal{P}| = |Q| |B| \,.
$$
But, clearly, $t \le \min \{ |Q|, |R|\} = |R|$ and $\tau \le \min \{ |A|, |B| \}$, thus we arrive to a contradiction in view of the assumption $|Q| \ge |A|$.
Finally, we need to consider the case when the third term in (\ref{tmp:03.12.2015_1}) dominates.
In the situation
$$
    t \tau |S_\tau| \ll |S_\tau| |R|
$$
%which is nonsense because of $t\le  \min \{ |Q|, |R|\}$.
and hence in view of (\ref{f:t_tau_range})
$$
    |R| |Q| |B| \ll |R|^2 \,.
$$
But this is a contradiction because $|Q| \ge |R|$ and $B$ is large enough.
This completes the proof of the lemma.
$\hfill\Box$

\bigskip

It is easy to see from the definition that $1\le d_* (A) \le |A|$.
The second inequality can be obtained if one put $Q=A$, $R=\{ 1\}$, $t=1$.

\bigskip

\begin{note}
In paper \cite{RR-NS}, Lemma 7 (see also \cite{SS1}, Lemma 27) the same result was obtained for the quantity
$$
    d(A) := \min_{C\neq \emptyset} \frac{|AC|^2}{|A||C|} \,.
$$
Clearly, $d_* (A) \le d(A)$. Indeed, just take $t=|C|$, $Q=AC$, and $R=C^{-1}$.
\end{note}

\begin{note}
    Let $A$ be a set and $\Pi = AA$ or $A/A$.
    By Katz--Koester inclusion \cite{kk} that is $|\Pi \cap \la \Pi| \ge |A|$ for any $\la \in A/A$ one
    has $d_* (\Pi) \le |\Pi|^3 / |A|^3$.
    The last estimate is usually better than ordinary $|\Pi \Pi|^2 / |\Pi|^2$ even if one applies  Pl\"{u}nnecke--Ruzsa inequality \cite{TV} (even for large subsets of $A$).
\end{note}

\bigskip

%There is
One can easily prove
an analog of Lemma \ref{l:d_r} in a dual form.
In the case for any sets $Q,R$ and a real number $t>0$ put
$$
    \Sym^+_{t} (Q,R) := \{ x ~:~ |Q\cap (x-R)| \ge t \}
$$
and consider the following quantity
\begin{equation}\label{f:d_r+}
    d_{+} (A) := \min_{t>0}\, \min_{\emptyset \neq Q,R \subset \R \setminus \{0\}
                    % ~:~ A\subseteq \Sym^\times_{t} (Q,R)
                    } \,
        \frac{|Q|^2 |R|^2}{|A| t^3} \,,
\end{equation}
    where the second minimum in (\ref{f:d_r+}) is taken over any $Q,R$ such that
    $A\subseteq \Sym^+_{t} (Q,R)$
    and $\max\{ |Q|,|R| \} \ge |A|$.
After that repeating the proof of Lemma \ref{l:d_r}, we need to estimate the cardinality of the set
$$
    S_\tau := \{ s \in  AB^{-1} ~:~ |A\cap sB| \ge \tau \} \,.
$$
So, we have arrived to the equation $ab^{-1} = s$, $s\in S_\tau$, $a\in A$, $b\in B$ and, further,
to the equation $q+r = sb$, $s\in S_\tau$, $b\in B$, $q\in Q$, $r\in R$.
It corresponds to the lines $l_{r,s} = \{ (x,y) ~:~ y+r=sx \}$ and Theorem \ref{t:SzT}, combining with the calculations of the rest of Lemma \ref{l:d_r} gives the result.

\bigskip

Thus, we have obtained an analog of Lemma \ref{l:d_r}.
%the result.

\begin{lemma}
    Let $A, B\subset \R$ be two finite sets.
    Then for any real number $\tau \ge 1$ one has
\begin{equation}\label{f:d_r+}
    \{ s \in AB^{-1} ~:~ |A\cap sB| \ge \tau \} \ll \frac{d_{+} (A) |A| |B|^2}{\tau^3} \,.
\end{equation}
\label{l:d_r+}
\end{lemma}

So, one can define a set $A \subset \R$ to be of (multiplicative) Szemer\'{e}di--Trotter type if inequality (\ref{f:d_r+}) holds for any $B \subset \R$ and every real number $\tau \ge 1$.

\bigskip

We will consider further generalizations of the quantities $d_* (A), d_{+} (A)$ in our forthcoming paper.

\bigskip
\section{The proof of the main results}
\bigskip

We need in two technical lemmas from \cite{KS_smd}.

Let $A\subset \R$, $0\not\in A$ be a finite set and $\tau>0$ be a real number.
%-л  вс кого $\la \in A/A$ определим set $A_\la = A\cap \la A \subseteq A$.
Let also $S'_\tau$ be a set
%+пределим
$$
 S_\tau'\subset S_\tau := \{ \la ~:~ \tau < |A_\la| \le 2 \tau \} \subseteq A/A
    %\,.
$$
and for any nonzero $\alpha_1,\alpha_2,\alpha_3$ and different $\la_1,\la_2,\la_3\in S_\tau'$ one has
$$
        \sigma (\a_1 A_{\la_1}, \a_2 A_{\la_2}, \a_3 A_{\la_3})\le\sigma \,.
$$

\begin{lemma}
  Let $A\subset \R$, $0\not\in A$ be a finite set, $\tau>0$ be a real number,
    \begin{equation}\label{cond:main_lemma}
        32\sigma \le \tau^2 \le |A+A| \sqrt{\sigma} \,,
    \end{equation}
    %a set
    and  $S_\tau'$,
    % and quantity
    $\sigma$ are defined above.
    Then
\begin{equation}\label{f:main_lemma}
    |A+A|^2 \ge
    %\frac{|S_\tau|}{2M} \left(\frac{\tau^2 M^2}{4} - \sigma M^4 \right)
     \frac{\tau^3 |S_\tau'|}{128 \sqrt{\sigma}}
     \,.
    %\,,
\end{equation}
%    where $\tilde{\sigma} \lesssim \sigma$.
\label{l:main_lemma}
\end{lemma}

\begin{lemma}
    Let $A\subset \R$, $0\not\in A$ be a finite set and $L\ge 1$ be a real number.
    Suppose that
%\begin{equation}\label{upA/A}
%    |A+A||A/A|\le|A|^3\,.
%%    |A/A| \le A^2/(\log^2|A|) \,,
%\end{equation}
%и
\begin{equation}\label{cond:Sol_energyA/A}
    |A+A|^2 |A/A| \le L |A|^4 \,.
\end{equation}
    Then there is  $\tau \ge \E^{\times} (A)/(2|A|^2)$ and some sets $S'_\tau \subseteq S_\tau \subseteq A/A$,
    $|S_\tau| \tau^2 \gtrsim \E^{\times} (A)$, $|S'_\tau|\ge|S_\tau|/2$ such that
    for any element $\la$ from $S'_\tau$ one has
%\begin{equation}\label{add_ener_low}
%\E^{+} (A_\la) \gtrsim \tau^3 L^{-4}
%\end{equation}
%and
\begin{equation}\label{A/A_low}
|A_\la/A_\la| \gtrsim \tau^2 L^{-16}\,.
\end{equation}

Similarly, if
\begin{equation}\label{cond:Sol_energyAA}
    |A+A|^2 |AA| \le L |A|^4
\end{equation}
then there exists $\tau \ge \E^{\times} (A)/(2|A|^2)$ and some sets $S'_\tau \subseteq S_\tau \subseteq A/A$,
    $|S_\tau| \tau^2 \gtrsim \E^{\times} (A)$, $|S'_\tau|\ge|S_\tau|/2$ such that
    for any $\la \in S'_\tau$,
    %, we have (\ref{add_ener_low}) and
    we have
\begin{equation}\label{AA_low}
|A_\la A_\la| \gtrsim \tau^2 L^{-16}\,.
\end{equation}
\label{l:smallL}
\end{lemma}

{\bf Proof of Theorem \ref{t:main2}.}
Consider the situation with $A/A$, because the case of $AA$ is similar.
By $\Pi$ denote $A/A$.
Without loosing of generality, suppose that $0\notin A$.
Now assume  that inequality (\ref{cond:Sol_energyA/A}) holds with some parameter $L$.
Let also  $|A/A|^3 \le L' |A|^{4}$.
Our task is to find a lower bound for the quantities $L$, $L'$.
Using Lemma \ref{l:smallL}, we find a number $\tau \ge \E^{\times} (A)/(2|A|^2)$ and a set
$S'_\tau \subseteq S_\tau \subseteq A/A$,
$|S_\tau| \tau^2 \gtrsim \E^{\times} (A)$, $|S'_\tau| \gtrsim |S_\tau|$ such that
for any element $\la$ from $S'_\tau$
one has
%SK $|A_\la / A_\la| \gtrsim L^{-16} \tau^2$.
$|A_\la / A_\la| \gtrsim \tau^2 L^{-16}$.
Using Katz--Koester inclusion, namely, $A_\la / A_\la \subseteq \Pi \cap \la \Pi^{-1}$, $\la \in \Pi$ (see \cite{kk}), we get
for any $\la \in S'_\tau$ that
\begin{equation*}\label{tmp:t}
%SK    |\Pi \cap \la \Pi^{-1}| \ge |A_\la / A_\la| \gtrsim L^{-16} \tau^2 := t \,.
    |\Pi \cap \la \Pi^{-1}| \ge |A_\la / A_\la| \gtrsim \tau^2  L^{-16}:= t \,.
\end{equation*}
In particular, $S'_\tau \subseteq \Sym^\times_t (\Pi,\Pi)$.
Because of  $S'_\tau \subseteq S_\tau$, we obtain
$$
    \sum_{a\in A} |A\cap a S'_\tau| = \sum_{\la \in S'_\tau} |A\cap \la A|
        \gg \tau |S'_\tau|% \gg \eta \tau |S_\tau|
$$
and hence there is $a\in A$ such that for the set  $A' := A\cap a S'_\tau$
one has
\begin{equation}\label{A'first_est}
    |A'| \gg \tau |S'_\tau| |A|^{-1}\,.
\end{equation}
We know that $S'_\tau \subseteq \Sym^\times_t (\Pi,\Pi)$.
Hence $A' \subseteq \Sym^\times_t (a\Pi,\Pi)$.
%Applying Lemma \ref{l:d_r}
Applying formula (\ref{f:d_r})
%SK with $Q=a\Pi$, $R=\Pi$ and $t=t$,
with $Q=a\Pi$, $R=\Pi$,
%It follows that
we obtain
\begin{equation}\label{tmp:d_*(A')}
%SK   d_* (A') \le \frac{|\Pi|^4 }{|A'| t^3} \ll \frac{|\Pi|^4 L^{48}}{|A'| \tau^6}
d_* (A') \lesssim \frac{|\Pi|^4 }{|A'| t^3} \ll \frac{|\Pi|^4 L^{48}}{|A'| \tau^6}
     \ll
        \frac{L^{48} |A| |\Pi|^4 }{|S_\tau| \tau^7} \,.
\end{equation}
Using Theorem \ref{t:s_SzT} and Lemma \ref{l:d_r}  as well as inequalities  (\ref{f:energy_KB}), (\ref{A'first_est}), (\ref{tmp:d_*(A')}), we get
$$
    |A+A| \ge |A'+A'| \gtrsim |A'|^{\frac{58}{37}} d_* (A')^{-\frac{21}{37}}
        \gtrsim
            (\tau |S_\tau| |A|^{-1})^{\frac{58}{37}}  (|S_\tau| \tau^7 L^{-48} |A|^{-1} |\Pi|^{-4})^{\frac{21}{37}}
$$
$$
    \gtrsim
        (\E^\times (A))^{\frac{79}{37}} \tau^{\frac{47}{37}} L^{-\frac{1008}{37}} |A|^{-\frac{79}{37}} |\Pi|^{-\frac{84}{37}}
            \gtrsim
                (\E^\times (A))^{\frac{126}{37}} L^{-\frac{1008}{37}} |A|^{-\frac{173}{37}} |\Pi|^{-\frac{84}{37}}
                    \ge
$$
$$
    \ge
        (|A|^4 / |\Pi|)^{\frac{126}{37}} L^{-\frac{1008}{37}} |A|^{-\frac{173}{37}} |\Pi|^{-\frac{84}{37}}
            \ge
$$
$$
            \ge
                |A|^{\frac{331}{37}} L^{-\frac{1008}{37}} |\Pi|^{-\frac{210}{37}}
                    \ge
                        L^{-\frac{1008}{37}} (L')^{-\frac{70}{37}} |A|^{\frac{51}{37}} \,.
$$
The last estimate is greater than $|A|^{4/3}$ by some power of $|A|$.
Easy calculations show that one can take any number less than
$\frac{5}{9813}$ for the constant $c$.
This concludes the proof.
$\hfill\Box$

\bigskip

{\bf Proof of Theorem \ref{t:main1}.}
Let $\Pi = A/A$, $|\Pi| = K|A|$.
%SK In the proof we can restrict ourselves considering just the case $|A|^{5/23} \lesssim K \lesssim |A|^{1/4}$.
In the proof we can restrict ourselves considering just the case $|A|^{5/23} \le K \le\gamma |A|^{1/4}$
where $\gamma>0$ is a small constant.
Without loosing of generality, suppose that $0\notin A$.
Using Dirichlet principle we find $\tau \ge |A|/(2K)$ 
%such that 
with 
$|S_\tau| \tau \gtrsim |A|^2$.
Consider two subsets $S'_\tau, S''_\tau$ of $S_\tau$ 
such that 
%with 
$|S'_\tau| = |S''_\tau| \ge |S_\tau|/2$ and for some parameter $\kappa \in (0,1]$ the following holds $|A_\la/A| \le \kappa |\Pi|$ for all $\la \in S'_\tau$ and
$|A_\la/A| \ge \kappa |\Pi|$ for any $\la \in S''_\tau$.
For any $\la \in S'_\tau$ one has
\begin{equation}\label{tmp:d_est-}
    d_* (A_\la) \le d(A_\la) \le \frac{\kappa^2 |\Pi|^2}{|A_\la| |A|} \le \kappa^2 |\Pi|^2 \tau^{-1} |A|^{-1} \,.
\end{equation}
Thus,
%using
applying
Corollary \ref{c:gen_sigma} and Lemma \ref{l:main_lemma}, we see that
$$
    |A+A|^2 \gg \tau^3 |S_\tau| (\kappa^2 |\Pi|^2 \tau^{-1} |A|^{-1})^{-1/6} \tau^{-5/6}
            =
        \tau^{7/3} |S_\tau| |A|^{1/6} |\Pi|^{-1/3} \kappa^{-1/3} \,,
$$
%SK provided by conditions (\ref{cond:main_lemma}) hold.
provided that conditions (\ref{cond:main_lemma}) hold.
Using inequalities $\tau \ge |A|/(2K)$ and $|S_\tau| \tau \gtrsim |A|^2$, we obtain
$$
    |A+A|^2 \gtrsim |A|^2 (|A|/K)^{4/3} |A|^{1/6} (|A|K)^{-1/3} \kappa^{-1/3}
        \gg
            |A|^{19/6} K^{-5/3} \kappa^{-1/3} \,.
$$
Hence
\begin{equation}\label{f:A+A_1}
    |A+A| \gtrsim |A|^{19/12} K^{-5/6} \kappa^{-1/6} \,.
\end{equation}

For the set $S''_\tau$ we use the arguments as in the proof of Theorem \ref{t:main2}.
Using Katz--Koester inclusion, namely, $A_\la/A \subseteq \Pi \cap \la \Pi^{-1}$, we get
for any $\la \in S''_\tau$ that
\begin{equation*}\label{tmp:t}
    |\Pi \cap \la \Pi^{-1}| \ge |A_\la /A| \ge \kappa |\Pi| := t \,.
\end{equation*}
In particular, $S''_\tau \subseteq \Sym^\times_t (\Pi,\Pi)$.
Because of  $S''_\tau \subseteq S_\tau$, we obtain
$$
    \sum_{a\in A} |A\cap a S''_\tau| = \sum_{\la \in S''_\tau} |A\cap \la A|
%SK        \gg \tau |S''_\tau|% \gg \eta \tau |S_\tau|
        \gg \tau |S''_\tau|\gg \eta \tau |S_\tau|
$$
and hence there is $a\in A$ such that for the set  $A' := A\cap a S''_\tau$
one has
\begin{equation}\label{A'first_est'}
%SK    |A'| \gg \tau |S''_\tau| |A|^{-1}\,.
    |A'| \gg \tau |S_\tau| |A|^{-1}\,.
\end{equation}
We know that $S''_\tau \subseteq \Sym^\times_t (\Pi,\Pi)$.
Hence $A' \subseteq \Sym^\times_t (a\Pi,\Pi)$.
Applying formula (\ref{f:d_r})
%SK with $Q=a\Pi$, $R=\Pi$ and $t=t$,
with $Q=a\Pi$, $R=\Pi$,
%It follows that
we obtain
\begin{equation}\label{tmp:d_*(A')'}
%SK    d_* (A') \le \frac{|\Pi|^4 }{|A'| t^3} \ll \frac{|\Pi|}{|A'| \kappa^3}
    d_* (A') \le \frac{|\Pi|^4 }{|A'| t^3} = \frac{|\Pi|}{|A'| \kappa^3}
     \ll
        \frac{|A| |\Pi| }{\kappa^3 |S_\tau| \tau} \,.
\end{equation}
Using Theorem \ref{t:s_SzT} and Lemma \ref{l:d_r} as well as inequalities  (\ref{f:energy_KB}), (\ref{A'first_est'}), (\ref{tmp:d_*(A')'}), we get
$$
    |A+A| \ge |A'+A'| \gtrsim |A'|^{\frac{58}{37}} d_* (A')^{-\frac{21}{37}}
        \gtrsim
            (\tau |S_\tau| |A|^{-1})^{\frac{58}{37}}  (\kappa^3 \tau |S_\tau| |A|^{-1} |\Pi|^{-1})^{\frac{21}{37}}
$$
\begin{equation}\label{f:A+A_2}
%SK    \ge
   =
        (|S_\tau| \tau)^{\frac{79}{37}} |A|^{-\frac{79}{37}} |\Pi|^{-\frac{21}{37}} \kappa^{\frac{63}{37}}
            \gtrsim
                |A|^{\frac{79}{37}} |\Pi|^{-\frac{21}{37}} \kappa^{\frac{63}{37}}
                    \ge
                        |A|^{\frac{58}{37}} K^{-\frac{21}{37}} \kappa^{\frac{63}{37}} \,.
\end{equation}
Combining bound (\ref{f:A+A_2}) with (\ref{f:A+A_1}), we find that the optimal choice of $\kappa$ is
\begin{equation}\label{f:kappa_def}
    \kappa = |A|^{\frac{7}{830}} K^{-\frac{59}{415}} \le 1
\end{equation}
%SK because in view of (\ref{f:A+A_1}) or just by
because $|A|^{5/23} \le K$.
%the first part
%SKinequality (\ref{f:main_1})
%SKof Theorem \ref{t:main}, we can trivially  assume that $K \ge |A|^{7/118}$.
%SK Substituting the last inequality into (\ref{f:A+A_1}), we obtain
%IS 
Substituting the last inequality into (\ref{f:A+A_1}), we obtain
$$
    |A+A| \gtrsim |A|^{19/12} K^{-5/6} (|A|^{\frac{7}{830}} K^{-\frac{59}{415}})^{-1/6}
        =
            |A|^{\frac{1313}{830}} K^{-\frac{336}{415}} \,.
$$

The only we need to check conditions (\ref{cond:main_lemma}).
The inequality $\tau^2 \ge 32 \sigma$ easily follows from (\ref{tmp:d_est-}) and inequality $K \lesssim |A|^{1/4}$.
Indeed by Corollary \ref{c:gen_sigma} and bounds (\ref{tmp:d_est-}), $\tau \ge |A|/(2K)$, we have
$$
%SK \sigma  \le (K^2 |A| \tau^{-1})^{1/3} \tau^{5/3} \ll \tau^2 \,.
 \sigma  \le (K^2 |A| \tau^{-1})^{1/3} \tau^{5/3} \ll \gamma^{2/3}\tau^2 \,,
$$
and $\sigma\le\tau^2/32$ if $\gamma$ is small enough.
It remains to check $\tau^2 \le |A+A| \sqrt{\sigma}$.
%If inequality (\ref{f:main_lemma}) then we are done, if not then
We have taken $\sigma = \tau^{4/3} K^{2/3} |A|^{1/3} \kappa^{2/3}$.
Thus we need to insure in the  inequality
\begin{equation}\label{tau8}
    \tau^8 \le |A+A|^6 K^2 |A| \kappa^2 \,.
\end{equation}
By
%inequality
bound
(\ref{f:Solymosi}) one has $|A+A|^2 \gg |A|^3 K^{-1} \log^{-1} |A|$.
In addition, in view of (\ref{f:kappa_def}) and because of $K \ll |A|^{1/4}$, we have
$\kappa \gg |A|^{-\frac{9}{332}} \ge |A|^{-\frac{3}{100}}$.
Thus, 
$$
|A+A|^6 K^2 |A| \kappa^2 \gg |A|^{10-3/50} K^{-1}\log^{-3} |A| \gg |A|^9, 
$$
and (\ref{tau8}) is true for large $|A|$ since $\tau\le|A|$.
This concludes the proof.
$\hfill\Box$

\bigskip
\section{Sum--products results with energies}
\label{sec:sp-energies}
\bigskip

In the section  we  prove sum--products results, which have deal just with the energies of sets but not with its sumsets or product sets.

%We begin
Let us start
with a lemma which can be interesting in its own right.

\begin{lemma}
    Let $A, P \subset \R$ be two sets.
    Put
$$
    \sigma_* := \sum_{x\in P} |A \cap x A| \,.
$$
    Then there is $A' \subseteq A$
    %and a number $q$
    such that
    $A'$ has SzT--type with $d_* (A') \lesssim \frac{|P|^2 |A|^2 |A'|^2}{\sigma_*^3}$
    and $|A'| \gtrsim \sigma_* |P|^{-1}$.
    Similarly, put
$$
    \sigma_+ := \sum_{x\in P} |A \cap (x + A)| \,.
$$
    Then there exists $A''\subseteq A$
    %and a number $q$
    such that
    $A''$ has SzT--type with $d_+ (A'') \lesssim \frac{|P|^2 |A|^2 |A''|^2}{\sigma_+^3}$
    and $|A''| \gtrsim \sigma_+ |P|^{-1}$.
\label{l:A',A''}
\end{lemma}
{\bf Proof.}
We have
$$
    \sigma_* = \sum_{x\in A} |P \cap x A^{-1}|
$$
and thus by the pigeonholing principle there is  a set $A'\subseteq A$ and a number $q \le |A|$ such that
$|A'| q \sim \sigma_*$ and $q< |P \cap x A^{-1}| \le 2 q$ for any $x\in A'$.
%differ at most twice on $A'$.
Because $q\le |P|$, we have $|A'| \gtrsim \sigma_* |P|^{-1}$.
Using Lemma \ref{l:d_r} with $Q=P$ and $R=A$  we see the set $A'$ has SzT--type such that $d_* (A')$ does not exceed
$$
    d_* (A') \ll \frac{|P|^2 |A|^2}{q^3 |A'|} \lesssim \frac{|P|^2 |A|^2 |A'|^2}{\sigma_*^3}
$$
as required.
By similar arguments and an application of Lemma \ref{l:d_r+} instead of Lemma \ref{l:d_r}, we obtain the existence of the set $A''$.
This completes the proof.
$\hfill\Box$

\bigskip

Now
%let us
we are ready to formulate the main result of the section, which shows that any set either has multiplicative energy or there is a large subset with small additive energy and visa versa.
Similar results were obtained in \cite{BW} but as we said in the introduction we do not use the Balog--Szemer\'{e}di--Gowers Theorem in the proof.

\begin{theorem}
    Let $A \subset \R$ be a set.
    Then there is $A_1 \subseteq A$ such that $|A_1| \gtrsim \E^\times (A) |A|^{-2}$ and
\begin{equation}\label{f:sp_E_E_1}
    \E^{+} (A_1) \E^\times (A) \lesssim |A_1|^{7/2} |A|^2 \,.
\end{equation}
    Similarly, there is $A_2 \subseteq A$ such that $|A_2| \gtrsim \E^+ (A) |A|^{-2}$ and
\begin{equation}\label{f:sp_E_E_2}
    \E^{\times} (A_2) \E^{+} (A) \lesssim |A_2|^{7/2} |A|^2 \,.
\end{equation}
\label{t:sp_E_E}
\end{theorem}

{\bf Proof.}
Put
$$
    \E^\times_3 (A) := \sum_{x} |A\cap x A|^3 \,.
$$
By the pigeonhole principle there is $P\subseteq A/A$ and a number $\D$ such that
$\D^3 |P| \sim \E^*_3 (A)$ and $\D < |A\cap x A| \le 2 \D$  for any $x\in P$.
Applying Lemma \ref{l:A',A''} with $\sigma_* \sim \D |P|$, we find a set $A_1\subseteq A$,
$|A_1| \gtrsim \D $ such that
$d_*(A_1) \lesssim \frac{|A|^2 |A_1|^2}{|P| \D^3}$.
We have $\D \gtrsim \E^\times_3 (A) \D^{-2} |P|^{-1}$ and hence by the
Cauchy--Schwarz inequality, we get $|A_1| \gtrsim \E^\times (A)^2 |A|^{-2} \D^{-2} |P|^{-1}$.
Next,
$$
    \E^\times (A) \ge \sum_{x\in P} |A\cap x A|^2 \ge \D^2 |P|\,.
$$
Therefore, $|A_1| \gtrsim \E^\times (A) |A|^{-2}$. Using Corollary \ref{c:gen_sigma}, we get
$$
    (\E^{+} (A_1))^2 \E^\times_3 (A) \lesssim (\E^{+} (A_1))^2 |P| \D^3 \ll |A_1|^7 |A|^2  \,.
$$
Finally,
%by
applying
the Cauchy--Schwarz inequality again, we obtain
$$
    \E^{+} (A_1) \E^\times (A) \lesssim |A_1|^{7/2} |A|^2
$$
as required.
By similar arguments we obtain the existence of the set $A_2$.
This completes the proof.
$\hfill\Box$

\bigskip

Now we can prove Theorem \ref{t:BW'} from the introduction.

\begin{cor}
    Let $A\subset \R$ be a finite set and $\d = 1/5$.
    Then there are two disjoint subsets $B$ and $C$ of $A$ such that $A = B\sqcup C$ and
$$
    \max\{ \E^{+} (B), \E^{\times} (C)\} \lesssim  |A|^{3-\d} \,.
$$
\end{cor}
{\bf Proof.}
    Let $M\ge 1$ be a parameter which we will choose later.
    Our arguments is a sort of an algorithm.
    We construct a decreasing sequence of sets $C_1=A \supseteq C_2 \supseteq \dots \supseteq C_k$ and an increasing sequence of sets $B_0 = \emptyset \subseteq B_1 \subseteq \dots \subseteq B_{k-1} \subseteq A$ such that
    for any $j=1,2,\dots, k$ the sets $C_j$ and $B_{j-1}$  are disjoint and moreover $A = C_j \sqcup B_{j-1}$.
    If at some step $j$ we have $\E^{\times} (C_j) \le |A|^3 / M$ then we stop our algorithm putting
    $C=C_j$, $B = B_{j-1}$, and $k=j-1$.
    In the opposite situation where $\E^{\times} (C_j) > |A|^3 / M$ we apply Theorem \ref{t:sp_E_E}
to the set $C_j$, finding the subset $D_j$ of $C_j$ such that $|D_j| \gtrsim |A|/M$ and
\begin{equation}\label{tmp:15.12.2015_1}
    \E^{+} (D_j) \lesssim |D_j|^{7/2} M |A|^{-1} \,.
\end{equation}
    After that we put $C_{j+1} = C_j \setminus D_j$, $B_j = B_{j-1} \sqcup D_j$ and repeat the procedure.
    Clearly, $B_k = \bigsqcup_{j=1}^k D_j$ and because of $|D_j| \gtrsim |A|/M$, we have $k\lesssim M$.
    Finally, by the H\"{o}lder inequality, Lemma \ref{1/4ineq} and (\ref{tmp:15.12.2015_1}), we get
$$
    (\E^{+} (B_k))^{1/4} \le \sum_{j=1}^k (\E^{+} (D_j))^{1/4} \lesssim (M |A|^{-1})^{1/4} \sum_{j=1}^k |D_j|^{7/8}
        \le
$$
$$
        \le
            (M |A|^{-1})^{1/4} \left( \sum_{j=1}^k |D_j| \right)^{7/8} k^{1/8}
                \lesssim
                    (M |A|^{-1})^{1/4} |A|^{7/8} M^{1/8}
                        =
                            M^{3/8} |A|^{5/8} \,.
$$
Hence
$$
    \E^{+} (B_k) \lesssim M^{3/2} |A|^{5/2} \,.
$$
Optimizing over $M$, that is choosing $M=|A|^{1/5}$, we obtain the result.
This completes the proof.
$\hfill\Box$

\bigskip

%    Using some results from \cite{s_mixed}, see section 7, one can improve lower bounds for the sizes of the sets $A_1$, %$A_2$ in Theorem \ref{t:sp_E_E}. Another way is to apply methods from \cite{S_BSzG}. We do not make such %calculations,just
%    Finally, let us show how to use the algorithm from the proof of the previous Corollary to increase the size of the %sets $A_1,A_2$.
%
%
We immediately get from Theorem~\ref{t:sp_E_E}  that
$$
\E^{+} (A_1)  \E^\times (A)\lesssim |A|^{11/2} \,
$$
and
$$
\E^{\times} (A_2) \E^{+} (A) \lesssim |A|^{11/2} \,.
$$
If $\E^\times (A)$ (respectively, $\E^{+} (A)$) is not too big, then it is not difficult to construct
larger sets $A_1$ and $A_2$ satisfying these inequalities.

\begin{cor}
    Let $A \subset \R$ be a set.
    Then there is $A_1 \subseteq A$ such that $|A_1| \gg (\E^\times (A))^{1/3}$ and
\begin{equation}\label{f:sp_E_E_1'}
    \E^{+} (A_1) \E^\times (A) \lesssim
        %|A_1|^{7/2} |A|^2 \,.
        |A|^{11/2} \,.
\end{equation}
    Similarly, there is $A_2 \subseteq A$ such that $|A_2| \gg (\E^+ (A))^{1/3}$ and
\begin{equation}\label{f:sp_E_E_2'}
    \E^{\times} (A_2) \E^{+} (A) \lesssim
        %|A_2|^{7/2} |A|^2 \,.
        |A|^{11/2} \,.
\end{equation}
%    In particular, $|A_1| \gg (\E^\times (A))^{1/3}$ and $|A_2| \gg (\E^+ (A))^{1/3}$.
\end{cor}
{\bf Proof.}
We proceed as in the proof of Theorem \ref{t:BW'}.

    We construct a decreasing sequence of sets $C_1=A \supseteq C_2 \supseteq \dots \supseteq C_k$ and an increasing sequence of sets $B_0 = \emptyset \subseteq B_1 \subseteq \dots \subseteq B_{k-1} \subseteq A$ such that
    for any $j=1,2,\dots, k$ the sets $C_j$ and $B_{j-1}$  are disjoint and moreover $A = C_j \sqcup B_{j-1}$.
    If at some step $j$ we have $|B_{j-1}| > (\E^{\times} (A))^{1/3} / 2$ then we stop our algorithm putting
    $A_1 = B_{j-1}$ and $k=j-1$.
    In the opposite situation where $|B_{j-1}| \le (\E^{\times} (A))^{1/3} / 2$
    we apply Theorem \ref{t:sp_E_E}
to the set $C_j$, finding the subset $D_j$ of $C_j$ such that $|D_j| \gtrsim \E^{\times}(C_j)/|C_j|^2$ and
$$
    \E^{+} (D_j) \lesssim |D_j|^{7/2}|C_j|^2 \E^{\times}(C_j)^{-1} \,.
$$
We observe, however,  that the inequality $\E^{\times}(B_{j-1})\le \E^{\times}(A)/8$ due to $|B_{j-1}| \le (\E^{\times} (A))^{1/3} / 2$ implies $\E^{\times}(C_j)\gg \E^{\times}(A)$. Therefore,
$$
|D_j| \gtrsim \E^{\times}(A)/|A|^2
$$
and
$$
    \E^{+} (D_j) \lesssim |D_j|^{7/2} |A|^2 \E^{\times}(A)^{-1} \,.
$$
Next we put $C_{j+1} = C_j \setminus D_j$, $B_j = B_{j-1} \sqcup D_j$ and repeat the procedure.

By Lemma \ref{1/4ineq}, we have
$$
    (\E^{+} (B_{k-1}))^{1/4} \le \sum_{j=1}^{k-1} (\E^{+} (D_j))^{1/4}
\lesssim |A|^{1/2}\E^{\times}(A)^{-1/4}\sum_{j=1}^{k-1} |D_j|^{7/8}
$$
$$
\lesssim |A|^{1/2}\E^{\times}(A)^{-1/4}\sum_{j=1}^{k-1}|D_j|
\left(\E^{\times}(A)/|A|^2\right)^{-1/8}
$$
$$
\le |A|^{1/2}\E^{\times}(A)^{-1/4}\E^{\times}(A)^{1/3}
\left(\E^{\times}(A)/|A|^2\right)^{-1/8}
$$
$$
=|A|^{3/4}\E^{\times}(A)^{-1/24} \le |A|^{11/8}\E^{\times}(A)^{-1/4}\,.
$$
So,
\begin{equation}
\label{B_k-1}
\E^{+} (B_{k-1}) \lesssim |A|^{11/2}\E^{\times}(A)^{-1}\,.
\end{equation}
Next,
\begin{equation}
\label{D_k}
    \E^{+} (D_k) \lesssim |D_k|^{7/2} |A|^2 \E^{\times}(A)^{-1}
\le |A|^{11/2} \E^{\times}(A)^{-1}\,.
\end{equation}
Since $A_1 = B_k = B_{k-1}\cup D_k$, we combine (\ref{B_k-1})
and (\ref{D_k}) to complete the proof of the first claim of
the corollary. The proof of the second claim is similar.
$\hfill\Box$

%\centerline{REFERENCES}
%\bigskip

{}
\bigskip
%\newpage

\noindent{S.V.~Konyagin\\
Steklov Mathematical Institute of Russian Academy of Sciences

{\tt konyagin@mi.ras.ru}

\bigskip

\noindent{I.D.~Shkredov\\
Steklov Mathematical Institute of Russian Academy of Sciences

{\tt ilya.shkredov@gmail.com}

\end{document}